\documentclass[12pt]{amsart}
\usepackage{amsmath,amscd}
\usepackage{xy}
\xyoption{all}
\newtheorem{thm}{Theorem}

\newtheorem{conjecture}[thm]{Conjecture}
\newtheorem*{thm*}{Theorem}
\newtheorem*{conjecture*}{Conjecture}

\newtheorem{prob}[thm]{Problem}
\numberwithin{equation}{section}
\numberwithin{thm}{section}

\begin{document}

\newcommand{\Hom}{\mathsf{Hom}}
\newcommand{\Aut}{\mathsf{Aut}}
\newcommand{\Inn}{\mathsf{Inn}}
\newcommand{\Ad}{\mathsf{Ad}}
\newcommand{\hol}{\mathsf{hol}}
\newcommand{\Uo}{\mathsf{U}(1)}
\newcommand{\Un}{\mathsf{U}(n)}
\newcommand{\SUn}{\mathsf{SU}(n)}
\newcommand{\PUn}{\mathsf{PU}(n)}
\newcommand{\B}{\mathbb{B}}
\newcommand{\tr}{\mathsf{tr}}
\newcommand{\Exp}{\mathsf{Exp}}
\newcommand{\Out}{\mathsf{Out}}
\newcommand{\Jac}{\mathsf{Jac}}
\newcommand{\Ker}{\mathsf{Ker}}
\newcommand{\Image}{\mathsf{Image}}
\newcommand{\Def}{\mathsf{Def}_{(G,G/H)}(\Sigma)}
\newcommand{\Ou}{\Out(\pi)}
\newcommand{\hpg}{\Hom(\pi,G)}
\newcommand{\cc}{\mathfrak{C}_\Sigma}
\newcommand{\hh}{\mathfrak{H}_\Sigma(G)}
\newcommand{\hpgg}{\Hom(\pi,G)/G}
\newcommand{\hpggmm}{\Hom(\pi,G)^{--}/G}
\newcommand{\Diff}{\mathsf{Diff}}
\newcommand{\Homeo}{\mathsf{Homeo}}
\newcommand{\mcg}{\mathsf{Mod}_{\Sigma}}
\newcommand{\Ht}{\mathsf{H}^2}
\newcommand{\Hthree}{\mathsf{H}^3}
\newcommand{\R}{\mathbb{R}}
\newcommand{\C}{\mathbb{C}}
\newcommand{\Z}{\mathbb{Z}}
\newcommand{\V}{\mathbb{V}}
\newcommand{\sut}{\mathsf{SU}(2)}
\newcommand{\sot}{\mathsf{SO}(2)}
\newcommand{\GL}{\mathsf{GL}}
\newcommand{\Sp}{\mathsf{Sp}}
\newcommand{\SL}{\mathsf{SL}}
\newcommand{\gltz}{\GL(2,\Z)}
\newcommand{\sltr}{\SL(2,\R)}
\newcommand{\sltc}{\SL(2,\C)}
\newcommand{\slthr}{\SL(3,\R)}
\newcommand{\rpt}{\R\mathsf{P}^2}
\newcommand{\psltr}{\mathsf{PSL}(2,\R)}
\newcommand{\pgltr}{\mathsf{PGL}(2,\R)}
\newcommand{\gltr}{\mathsf{GL}(2,\R)}
\newcommand{\psltc}{\mathsf{PSL}(2,\C)}
\newcommand{\cpo}{\C\mathsf{P}^1}
\newcommand{\Teich}{\mathfrak{T}_{\Sigma}}
\newcommand{\Is}{\mathcal{I}_\Sigma}
\newcommand{\g}{\mathfrak{g}}
\newcommand{\ggg}{\g_{\Ad\rho}}
\newcommand{\Diffo}{\mathsf{Diff}^0(\Sigma)}
\newcommand{\Homeoo}{\mathsf{Homeo}^0(\Sigma)}
\newcommand{\qf}{\mathcal{Q}_\Sigma}
\newcommand{\Sym}{\mathsf{Sym}}
\newcommand{\SSS}{\mathfrak{S}}
\newcommand{\Uu}{\mathbb{U}}
\newcommand{\Mod}{\mathfrak{M}_\Sigma}
\newcommand{\FF}{\mathfrak{F}}
\newcommand{\UU}{\mathfrak{U}}
\newcommand{\ZZ}{\mathfrak{Z}}
\newcommand{\HH}{\mathfrak{H}}
\newcommand{\dd}[1]{\partial_{#1}}

 \setcounter{tocdepth}{2}

\title[Surface Group Representations]
{Mapping Class Group Dynamics on Surface Group Representations}
\author{William M. Goldman} 
\address{Department of Mathematics \\ 
 University of Maryland \\ 
College Park, MD 20742}
\thanks{Goldman supported in part by NSF grants DMS-0103889 and
DMS-0405605} 
\email{wmg@math.umd.edu}
\date{\today}


\begin{abstract}
Deformation spaces $\hpgg$ of representations of the fundamental group
$\pi$ of a surface $\Sigma$ in a Lie group $G$ admit natural actions
of the mapping class group $\mcg$, preserving a Poisson
structure. When $G$ is compact, the actions are ergodic. In contrast
if $G$ is noncompact semisimple, the associated deformation space
contains open subsets containing the Fricke-Teichm\"uller space 
upon which $\mcg$ acts properly. Properness of the $\mcg$-action
relates to (possibly singular) locally homogeneous geometric structures
on $\Sigma$. We summarize known results and state open questions about
these actions.
\end{abstract}

\maketitle
\tableofcontents
\section*{Introduction}

A natural object associated to a topological surface
$\Sigma$ is the {\em deformation space\/} of representations of its
fundamental group $\pi=\pi_1(\Sigma)$ in a Lie group $G$. 
These spaces admit natural actions of the mapping class group $\mcg$
of $\Sigma$, and therefore determine linear representations of $\mcg$.

The purpose of this paper is to survey recent results on the dynamics
of these actions, and speculate on future directions in this subject.

The prototypes of this theory are two of the most basic spaces in Riemann 
surface theory: the Jacobian and the Fricke-Teichm\"uller space.
The {\em Jacobian\/} $\Jac(M)$ of a Riemann surface $M$ homeomorphic to 
$\Sigma$ identifies with the deformation space $\hpgg$ 
when $G$ is the circle $\Uo$.
The Jacobian parametrizes topologically trivial holomorphic complex
line bundles over $M$, but its topological type (and symplectic structure)
are invariants of the underlying topological surface $\Sigma$. The action
of $\mcg$ is the action of the integral symplectic group $\Sp(2g,\Z)$
on the torus $\R^{2g}/\Z^{2g}$, which is a measure-preserving chaotic (ergodic)
action.

In contrast, the {\em Teichm\"uller space\/} $\Teich$ ({\em Fricke
space\/} if $\partial\Sigma\neq\emptyset$) is comprised of equivalence
classes of {\em marked conformal structures\/} on $\Sigma$.  A marked
conformal structure is a pair $(M,f)$  where $f$ is a homeomorphism 
and $M$ is a Riemann surface.  Marked conformal structures
$(f_1,M_1)$ and $(f_2,M_2)$  are {\em equivalent\/} if there
is a biholomorphism $M_1\xrightarrow{h} M_2$ such that $h\circ f_1$ is
homotopic to $f_2$.  Denote the equivalence class of a marked
conformal structure $(f,M)$ by 
\begin{equation*}
\langle f,M\rangle\in\Teich.  
\end{equation*}
A marking $f$ determines a representation
of the fundamental group:
\begin{equation*}
\pi = \pi_1(\Sigma) \xrightarrow{f_*} \pi_1(M) \subset \Aut(\tilde M). 
\end{equation*}

By the uniformization theorem (at least when
$\chi(\Sigma)<0$), these identify with {\em marked hyperbolic
structures\/} on $\Sigma$, which in turn identify with conjugacy
classes of {\em discrete embeddings\/} of the fundamental group $\pi$
in the group $G=\pgltr$ of isometries of the hyperbolic plane. These classes 
form a connected component of $\hpgg$, which is homeomorphic to a cell
of dimension $-3\chi(\Sigma)$~\cite{doctoralthesis}.  
The mapping class group $\mcg$ acts properly on $\Teich$. 
The quotient orbifold
\begin{equation*}
\Mod := \Teich/\mcg 
\end{equation*}
is the {\em Riemann moduli space,\/} consisting of biholomorphism classes
of (unmarked) conformal structures on $\Sigma$.

Summarizing:
\begin{itemize}
\item When $G$ is compact, $\hpgg$ has 
nontrivial homotopy type, 
and the action of the mapping class group
exhibits nontrivial dynamics;
\item When $G=\pgltr$ (or more generally a noncompact semisimple Lie group),
$\hpgg$ contains open sets (like Teichm\"uller space)  which 
are contractible and admit a proper $\mcg$-action. Often these
open sets correspond to locally homogeneous geometric structures 
uniformizing $\Sigma$.
\end{itemize}
Thus dynamically complicated mapping class group actions accompany
nontrivial homotopy type of the deformation space. In general the dynamics exhibits
properties of these two extreme cases, as will be described in this
paper. 

\subsection*{Acknowledgments}
This paper is an expanded version of a lecture presented at the
Special Session ``Dynamics of Mapping Class Group Actions'' at the
Annual Meeting of the American Mathematical Society, January 6-11,
2005, in Atlanta, Georgia.  I am grateful to Richard Brown for
organizing this workshop, and the opportunity to lecture on this
subject. I am also grateful to Benson Farb for encouraging me to write
this paper, and to J\o rgen Andersen, David Dumas, Lisa Jeffrey, Misha
Kapovich, Fran\c cois Labourie, Dan Margalit, Howard Masur, Walter
Neumann, Juan Souto, Pete Storm, Ser Tan, Richard Wentworth, Anna
Wienhard and Eugene Xia for several suggestions and helpful comments.
I wish to thank the referee for a careful reading of the paper and
many useful suggestions.

\section{Generalities}
Let $\pi$ be a finitely generated group and $G$ a real algebraic
Lie group.  The set $\hpg$ of homomorphisms $\pi\longrightarrow G$ has
the natural structure of an affine algebraic set.  The group
\begin{equation*}
\Aut(\pi)\times\Aut(G) 
\end{equation*}
acts on $\hpg$ by left- and right- composition, 
preserving the algebraic structure: if $\alpha\in\Aut(\pi)$ and
$h\in\Aut(G)$ are automorphisms, 
then the action of $(\alpha,h)$ on $\rho\in\hpg$ is
the composition $h\circ \rho\circ\alpha^{-1}$:
\begin{equation*}
\pi \xrightarrow{\alpha^{-1}} \pi \xrightarrow{\rho} G \xrightarrow{h} G 
\end{equation*}

The {\em deformation space\/} is the quotient space of $\hpg$
(with the classical topology) by the subgroup $\Inn(G)$ 
of {\em inner automorphisms\/} of $G$, and is denoted $\hpgg$.
The action of the inner automorphism 
$\iota_\gamma$ determined by an element $\gamma\in\pi$ 
equals $\iota_{\rho(\gamma^{-1})}(\rho)$.
Therefore $\Inn(\pi)$ acts trivially on $\hpgg$ and
the induced action of $\Aut(\pi)$
on $\hpgg$ factors through the quotient
\begin{equation*}
\Ou := \Aut(\pi)/\Inn(\pi). 
\end{equation*}
When $\Sigma$ is a closed orientable surface with $\chi(\Sigma)<0$, then
the natural homomorphism
\begin{equation*}
\pi_0(\Diff(\Sigma)) \longrightarrow \Ou 
\end{equation*}
is an isomorphism.
The {\em mapping class group\/} $\mcg$ is the subgroup of $\Ou$
corresponding to {\em orientation-preserving\/} diffeomorphisms of $\Sigma$.

When
$\Sigma$ has nonempty boundary with components
$\partial_i\Sigma$, this deformation space admits 
a {\em boundary restriction map\/}
\begin{equation}\label{eq:relativecv}
\Hom(\pi_1(\Sigma),G)/G \longrightarrow 
\prod_{i\in\pi_0(\partial\Sigma)}
\Hom\big(\pi_1(\partial_i\Sigma)\big),G)/G. 
\end{equation}
The fibers of the boundary restriction map are the {\em relative character
varieties.\/} 
This action of $\mcg$ preserves this map.


\subsection{The Symplectic Structure}
These spaces possess {\em algebraic symplectic structures,\/}
invariant under $\mcg$.
For the moment we focus on the {\em smooth part\/} of
$\hpg$, which we define as follows. 
When $G$ is reductive, the subset $\hpg^{--}$
consisting of representations whose image does not lie in
a parabolic subgroup of $G$ is a smooth submanifold
upon which $\Inn(G)$ acts properly and freely. 
The quotient $\hpggmm$ is then a smooth  manifold,
with a $\mcg$-invariant symplectic structure.

The symplectic structure depends on a choice of a nondegenerate 
$\Ad$-invariant symmetric bilinear form $\B$ on the 
Lie algebra $\g$ of $G$ and an orientation on $\Sigma$.
The composition
\begin{equation*}
\pi \stackrel{\rho}\longrightarrow G  \stackrel{\Ad}\longrightarrow 
\Aut(\g)
\end{equation*}
defines a $\pi$-module $\ggg$.
The Zariski tangent space to $\hpg$ at a representation $\rho$ 
is the space $Z^1(\pi,\ggg)$  of {\em 1-cocycles.\/} 
The tangent space to the orbit $G\rho$ equals the subspace
$B^1(\pi,\ggg)$  of {\em 1-coboundaries.\/}
These facts are due to Weil~\cite{Weil}, 
see also Raghunathan~\cite{Raghunathan}.
If  $G$ acts properly and freely on a neighborhood of $\rho$ in $\hpg$,
then $\hpgg$ is a manifold near $[\rho]$ with tangent space 
$H^1(\pi,\ggg)$.
In that case a nondegenerate symmetric
$\Ad(G)$-invariant bilinear form 
\begin{equation*}
\g \times\g \xrightarrow{\B}\R 
\end{equation*}
defines a pairing of $\pi$-modules 
\begin{equation*}
\ggg \times\ggg \xrightarrow{\B}\R. 
\end{equation*}
Cup product using $\B$ as coefficient pairing 
defines a nondegenerate skew-symmetric pairing
\begin{equation*}
H^1(\pi,\ggg) \times
H^1(\pi,\ggg) \xrightarrow{\B_*(\cup)} H^2(\pi,\R) \cong \R
\end{equation*}
on each tangent space 
\begin{equation*}
T_{[\rho]}\hpgg \cong H^1(\pi,\ggg). 
\end{equation*}
Here the isomorphism $H^2(\pi,\R)\cong\R$ arises from the orientation
on $\Sigma$.  The resulting exterior 2-form $\omega_\B$ is
closed~\cite{nature}, and defines a symplectic structure on the smooth
part $\hpggmm$ of $\hpgg$. This topological definition makes it
apparent that $\omega_\B$ is $\mcg$-invariant.  In particular the
action preserves the measure $\mu$ defined by $\omega_\B$.  When $G$
is compact, the total measure is finite
(Jeffrey-Weitsman~\cite{JW1,JW2}, Huebschmann~\cite{Huebschmann}).

\subsection{The Complex  Case}

When $G$ is a complex Lie group, $\Hom(\pi,G)$ has a complex
algebraic structure preserved by the $\Aut(\pi)\times\Aut(G)$-action.
When $G$ is a complex semisimple Lie group, the above
construction, applied to a nondegenerate  $\Ad$-invariant {\em complex-bilinear
form\/} 
\begin{equation*}
\g \times\g \xrightarrow{\B}\C,  
\end{equation*}
determines a {\em complex-symplectic structure\/} on $\hpggmm$, 
that is, a closed
nondegenerate holomorphic $(2,0)$-form.  This complex-symplectic
structure is evidently $\mcg$-invariant.  
For a discussion of this structure when $G = \sltc$,
see \cite{csymp}.

The choice of a marked conformal structure on $\Sigma$ 
determines a {\em hyper-K\"ahler structure\/} on
$\hpgg$ subordinate to this complex-symplectic structure.

A complex-symplectic structure on a $4m$-dimensional real manifold $V$
is given by an integrable
almost complex structure $J$ and a closed nondegenerate skew-symmetric
bilinear form 
\begin{equation*}
TM \times TM\xrightarrow{\Omega}\C 
\end{equation*}
which is complex-bilinear with respect to $J$. 
Alternatively, it is defined by a reduction of the structure
group of the tangent bundle $TV$ from $\GL(4m,\R)$
to the  subgroup 
\begin{equation*}
\Sp(2m,\C) \subset \GL(4m,\R).
\end{equation*}

A hyper-K\"ahler structure further reduces the structure group 
of the tangent bundle from $\Sp(2m,\C)$ to its maximal compact subgroup
$\Sp(2m)\subset\Sp(2m,\C)$. All of these structures are required to satisfy
certain integrability conditions.
A hyper-K\"ahler structure  subordinate to a complex-symplectic
structure $(\Omega,J)$ is defined
by a Riemannian metric $g$ and integrable almost complex structures $I,K$
such that:
\begin{itemize}
\item  $g$ is K\"ahlerian with respect to each of $I,J,K$,
\item the complex structures $I,J,K$ satisfy the quaternion identities, 
\item $\Omega(X,Y) = - g(I X,Y) + i\, g(K X)$
for $X,Y\in TM$.
\end{itemize}
Goldman-Xia~\cite{GoldmanXia}, \S 5 describes this structure in detail
when $G = \GL(1,\C)$.

From this we can associate to every point in Teichm\"uller space $\Teich$  
a compatible hyper-K\"ahler structure on the complex-symplectic
space $\hpggmm$. 
However the hyper-K\"ahler structures are not $\mcg$-invariant.

\subsection{Singularities of the deformation space}
In general the spaces $\hpg$ and $\hpgg$ are not manifolds, 
but their local structure admits a very explicit cohomological
description. For convenience assume that $G$ is reductive algebraic
and that $\rho$ is a {\em reductive representation,\/} that is,
its image $\rho(\pi)$ is Zariski dense in a reductive subgroup of $G$.
For $\rho\in\hpg$, denote the centralizer of $\rho(\pi)$ by $\ZZ(\rho)$
and the center of $G$ by $\ZZ$.

A representation $\rho\in\hpg$ is a {\em singular point \/} of $\hpg$
if and only if 
\begin{equation*}
\dim(\ZZ(\rho)/\ZZ) > 0.   
\end{equation*}
Equivalently, the isotropy group of $\Inn(G)$ at $\rho$
is not discrete, that is, the action of $\Inn(G)$
at $\rho$ is not {\em locally free.\/}

The Zariski tangent space $T_{\rho}\hpg$ equals the space $Z^1(\pi;\ggg)$
of $\ggg$-valued 1-cocycles on $\pi$.
The tangent space  
to the orbit  $G\cdot\rho$ 
equals the
subspace $B^1(\pi;\ggg)$ of coboundaries. 
Thus the {\em Zariski normal space\/}
at $\rho$ to the orbit $G\cdot\rho$ in $\hpg$ equals the cohomology group
$H^1(\pi;\ggg)$. 

Here is a heuristic interpretation. Consider an analytic path
$\rho_t\in\hpg$ with $\rho_0 = \rho$.
Expand it as a power series in $t$:
\begin{equation}\label{eq:powerseries}
\rho_t(x) = \exp\big( u_0(x) t + u_2(x) t^2 + u_3(x) t^3 + \dots\big) \rho(x) 
\end{equation}
where 
\begin{equation*}
\pi\xrightarrow{u_n} \g
\end{equation*}
for $n\ge 0$. The condition 
\begin{equation}\label{eq:homo}
\rho_t(xy) = \rho_t(x) \rho_t(y) 
\end{equation}
implies that the tangent vector $u=u_0$ satisfies 
the {\em cocycle condition\/}
\begin{equation}\label{eq:cocyclecondition}
u(xy) = u(x) + \Ad\rho(x) u(y),
\end{equation}
(the {\em linearization \/} of \eqref{eq:homo}. 
The vector space of solutions of \eqref{eq:cocyclecondition}
is the space $Z^1(\pi;\ggg)$ 
of $\ggg$-valued {\em 1-cocycles\/} of $\pi$.

The Zariski tangent space to the orbit $G\cdot\rho$ equals the
subspace $B^1(\pi,\ggg)\subset Z^1(\pi,\ggg)$ consisting of {\em
1-coboundaries.\/}
Suppose that a path $\rho_t$ in $\hpg$ is induced by
a conjugation by a path $g_t$
\begin{equation*}
\rho_t(x) = g_t \rho(x) g_t^{-1},
\end{equation*}
where $g_t$ admits a power series expansion
\begin{equation*}
g_t = \exp (v_1 t + v_2 t^2 + \dots), 
\end{equation*}
where $v_1,v_2,\dots\in \g$.
Thus the tangent vector to $\rho_t$ is tangent to the orbit $G\cdot\rho$.
Expanding the power series, this tangent vector equals
\begin{equation*}
u(x) = v_1 - \Ad\rho(x) v_1,
\end{equation*}
that is, $u= \delta v_1 \in B^1(\pi;\ggg)$ is a coboundary.

Let $u \in T_\rho\hpg = Z^1(\pi;\ggg)$ 
be a tangent vector to $\hpg$ at $\rho$.
We give necessary and sufficient conditions that $u$ be tangent to an 
analytic path of representations.

Solving the equation \eqref{eq:homo} to second order gives:
\begin{equation}\label{eq:secondorder}
u_2(x) - u_2(x y) + \Ad\rho(x) u_2(y) = \frac12 [ u(x), \Ad\rho(x) u(y)].
\end{equation}
Namely, the function, 
\begin{align}\label{eq:cup}
\pi \times \pi & \longrightarrow \g \notag\\
(x,y) & \longmapsto \frac12 [ u(x), \Ad\rho(x) u(y)]
\end{align}
is a $\ggg$-valued 2-cochain on $\pi$,
This 2-cochain is the coboundary $\delta u_2$
of the 1-cochain $\pi\xrightarrow{u_2}\ggg$.
Similarly there are conditions on the coboundary of $u_n$ in terms
of the lower terms in the power series expansion \eqref{eq:powerseries}.

The operation \eqref{eq:cup} has a cohomological interpretation as follows.
$\pi$ acts on $\g$ by Lie algebra automorphisms, 
so that Lie bracket defines a pairing of $\pi$-modules
\begin{equation*}
\ggg\times\ggg\xrightarrow{[,]}\ggg. 
\end{equation*}
The Lie algebra of $\ZZ(\rho)$ equals $H^0(\pi;\ggg)$.
The linearization of the action of $\ZZ(\rho)$ is given by
the cup product on $H^1(\pi;\ggg)$ with $[,]$  as coefficient
pairing:
\begin{equation*}
H^0(\pi;\ggg) \times H^1(\pi;\ggg)
\stackrel{[,]_*(\cup)}\longrightarrow H^1(\pi;\ggg). 
\end{equation*}

Now consider the cup product of $1$-dimensional classes.
The bilinear form
\begin{equation*}
H^1(\pi;\ggg) \times H^1(\pi;\ggg)
\stackrel{[,]_*(\cup)}\longrightarrow H^2(\pi;\ggg). 
\end{equation*}
is symmetric; let $Q_\rho$ be the associated quadratic form.

Suppose $u$ is tangent to an analytic path. Solving
\eqref{eq:powerseries} to second order (as in \eqref{eq:secondorder}
and \eqref{eq:cup}) implies that
\begin{equation*}
[,]_*(\cup)([u],[u]]) = \delta u_2, 
\end{equation*}
that is,
\begin{equation}\label{eq:necessarycondition}
Q_\rho([u])=0. 
\end{equation}

Under the above hypotheses, the necessary condition
\eqref{eq:necessarycondition} is also sufficient. In fact,
by Goldman-Millson~\cite{defth}, $\rho$ has a neighborhood $N$ in 
$\hpg$ analytically equivalent to a neighborhood
of $0$ of the cone $C_\rho$ in $Z^1(\pi;\ggg)$ defined by the homogeneous
quadratic function
\begin{align*}
Z^1(\pi;\ggg) & \longrightarrow H^2(\pi;\ggg) \\
 u & \longmapsto Q_\rho([u]).
\end{align*}
Then the germ of $\hpgg$ at $[\rho]$ is the quotient of this cone by
the isotropy group $\ZZ(\rho)$.
(These spaces are special cases of {\em symplectic stratified spaces\/}
of Sjamaar-Lerman~\cite{SjamaarLerman}.)

An explicit {\em exponential mapping\/}
\begin{equation*}
N \xrightarrow{~\Exp_\rho~} \hpg
\end{equation*}
was constructed by Goldman-Millson~\cite{defkah} using the Green's operator 
of a Riemann surface $M$ homeomorphic to $\Sigma$.

The subtlety of these constructions is underscored by the following
false argument, which seemingly proves that the Torelli subgroup
of $\mcg$ acts identically on the whole component of $\hpgg$ containing
the trivial representation. This is easily seen to be false, for
$G$ semisimple.

Here is the fallacious argument. 
The trivial representation $\rho_0$ is fixed by all of $\mcg$.
Thus $\mcg$ acts on the analytic germ of $\hpgg$ at $\rho_0$.
At $\rho_0$, the coefficient module $\ggg$ is 
trivial, and the tangent space corresponds to ordinary (untwisted)
cohomology:
\begin{equation*}
T_{\rho_0}\hpg = Z^1(\pi;\g)  = Z^1(\pi)\otimes \g. 
\end{equation*}
The quadratic form is just the usual cup-product pairing, so
any homologically trivial automorphism $\phi$ fixes the quadratic cone $N$
pointwise. By Goldman-Millson~\cite{defth}, the analytic
germ of $\hpg$ at $\rho_0$ is equivalent to the quadratic cone $N$.
{\em Therefore $[\phi]$ acts trivially on an open neighborhood
of $\rho$ in $\hpg$.\/} By analytic continuation, $[\phi]$ acts trivially
on the whole component of $\hpg$ containing $\rho$.

The fallacy arises 
because the identification $\Exp_\rho$ of a neighborhood $N$
in the quadratic cone with the germ of $\hpg$ at $\rho$ depends on a
choice of Riemann surface $M$. 
Each point $\langle f,M\rangle\in \Teich$ determines
an exponential map 
$\Exp_{\rho,\langle f,M\rangle}$ 
from the germ of the quadratic cone to
$\hpg$, and these are {\em not \/} invariant under $\mcg$.
In particular, no family of isomorphisms of the analytic germ of
$\hpg$ at $\rho_0$ with the quadratic cone $N$ is $\mcg$-invariant.

\begin{prob} Investigate the dependence 
of $\Exp_{\rho,\langle f,M\rangle}$ 
on the marked Riemann surface $\langle f,M\rangle$.
\end{prob}

\subsection{Surfaces with boundary}
When $\Sigma$ has nonempty boundary, an $\Ad$-invariant 
inner product $\B$ on $\g$ 
and an orientation on $\Sigma$
determines a {\em Poisson structure\/}
(Fock-Rosly~\cite{FockRosly},
Guruprasad-Huebschmann-Jeffrey-Weinstein~\cite
{GuruprasadHuebschmannJeffreyWeinstein}). 
The {\em symplectic leaves\/}
of this Poisson structure are the level sets of the boundary 
restriction map \eqref{eq:relativecv}.

For each component $\partial_i\Sigma$ of $\partial\Sigma$, 
fix a conjugacy class $C_i\subset G$. 
The subspace 
\begin{equation}
\hpgg_{(C_1,\dots,C_b)} \subset \hpgg 
\end{equation}
consisting of $[\rho]$ such that
\begin{equation}\label{eq:relcv}
\rho(\partial_i\Sigma)\subset C_i 
\end{equation}
has a symplectic structure.
(To simplify the discussion we assume that it is a smooth submanifold.)
De Rham cohomology with twisted coefficients in $\ggg$ is 
naturally isomorphic with group cohomology of $\pi$. 
In terms of De Rham cohomology, the tangent space at $[\rho]$
to $\hpgg_{(C_1,\dots,C_b)}$
identifies with 
\begin{align*}
\Ker\bigg( H^1(\Sigma;& \ggg) \to   H^1(\partial\Sigma;\ggg)\bigg) \\
& \cong \; 
 \Image\bigg( H^1(\Sigma,\partial\Sigma;\ggg)
\to H^1(\Sigma;\ggg) \bigg).
\end{align*}
The cup product pairing
\begin{equation*}
H^1(\Sigma;\ggg) \times H^1(\Sigma,\partial\Sigma;\ggg) 
\xrightarrow{\B_*(\cup)} H^2(\Sigma,\partial\Sigma;\R) 
\end{equation*}
induces a symplectic structure on 
$\hpgg_{(C_1,\dots,C_b)}$.

Given a (possibly singular) foliation $\FF$ of a manifold $X$ 
by symplectic manifolds, the Poisson structure is defined as follows. 
For functions $f,g\in C^\infty(X)$, their Poisson bracket is a function
$\{f,g\}$ on $X$ defined as follows.
Let $x\in X$ and let $L_x$ be the leaf of $\FF$ containing $x$.
Define the value of $\{f,g\}$ at $x$ as the Poisson bracket 
\begin{equation*}
\{f|_{L_x},g|_{L_x}\}_{L_x}, 
\end{equation*}
where $\{,\}_{L_x}$ denotes the Poisson bracket operation on the
symplectic manifold $L_x$, and 
$f|_{L_x},g|_{L_x}\in C^\infty(L_x)$,
are the restrictions of $f,g$ to $L_x$.

The examples below exhibit {\em exterior bivector fields\/} $\xi$
representing the Poisson structure. 
If  $f,g\in C^\infty(X)$, their Poisson bracket $\{f,g\}$
is expressed as an interior product of $\xi$
with the exterior derivatives of $f,g$:
\begin{equation*}
\{ f, g\} = \xi \cdot (df \otimes dg).
\end{equation*}
In local coordinates $(x^1,\dots,x^n)$, write 
\begin{equation*}
\xi = \sum_{i,j} \xi^{i,j} \frac{\partial}{\partial x_i} 
\wedge \frac{\partial}{\partial x_j} 
\end{equation*}
with $\xi^{i,j}=-\xi^{j,i}$. Then
\begin{align*}
\{ f, g\} & = \sum_{i,j}  \bigg(\xi^{i,j} 
\frac{\partial}{\partial x_i}\wedge
\frac{\partial}{\partial x_j}\bigg) 
\cdot \bigg(\frac{\partial f}{\partial x_i} dx^i \otimes
\frac{\partial g}{\partial x_j} dx^j\bigg) \\
& = \sum_{i,j}  \xi^{i,j}
\bigg( 
\frac{\partial f}{\partial x_i} \frac{\partial g}{\partial x_j} - 
\frac{\partial f}{\partial x_j} \frac{\partial g}{\partial x_i}\bigg)
\end{align*}
\subsection{Examples of relative $\sltc$-character varieties}
We give a few explicit examples, when $G=\sltc$, and $\Sigma$ is a
three-holed or four-holed sphere, or a one-holed or two-holed torus.
Since generic conjugacy classes in $\sltc$ are determined by the trace function
\begin{equation*}
\sltc \xrightarrow{\tr} \C 
\end{equation*}
the relative character varieties are level sets of the mapping
\begin{align*}
\hpgg & \longrightarrow \C^b   \\
[\rho] & \longmapsto \bigg[
\tr\big(\rho(\partial_i(\Sigma)\big)\bigg]_{i=1,\dots,b}
\end{align*}

\subsubsection{The three-holed sphere}

When $\Sigma$ is a three-holed sphere, 
its fundamental group admits a presentation
\begin{equation*}
\pi = \langle A, B, C \mid  ABC = 1 \rangle 
\end{equation*}
where $A, B, C$ correspond to the three components of $\partial
\Sigma$.  Here is the fundamental result for $\sltc$-character varieties of
a rank two free group:

\begin{thm*}[Vogt~\cite{Vogt}-Fricke~\cite{FrickeKlein}]  
The map
\begin{align*}
\hpgg & \longrightarrow \C^3 \\
[\rho] &\longmapsto 
\bmatrix
\tr\big(\rho(A)\big) \\ \tr\big(\rho(B)\big) \\ \tr\big(\rho(C)\big) 
\endbmatrix
\end{align*}
is an isomorphism of affine varieties. 
\end{thm*}
\noindent
In particular, the symplectic leaves are just points.
See \cite{fricke} for an elementary proof..

\subsubsection{The one-holed torus}\label{subsec:oneholedtorus}
When $\Sigma$ is a one-holed torus, its fundamental group admits a presentation
\begin{equation*}
\pi = \langle X, Y, Z, K  \mid  X Y Z = 1, K = X Y X^{-1} Y^{-1} \rangle 
\end{equation*}
where $X,Y$ are simple loops intersecting once, 
and $K = X Y X^{-1} Y^{-1}$ corresponds to the boundary.
Presenting the interior of $\Sigma$ as the quotient
\begin{equation*}
\mathsf{int}\big(\Sigma\big) = \big(\R^2-\Z^2\big)/\Z^2 
\end{equation*}
the curves $X, Y$ correspond to the $(1,0)$ and $(0,1)$-curves
respectively. Once again the Vogt-Fricke theorem implies that
$\hpgg\cong\C^3$, with coordinates
\begin{align*}
x & = \tr\big( \rho(X)\big)  \\
y & = \tr\big( \rho(Y)\big) \\ 
z & = \tr\big( \rho(Z)\big).  
\end{align*}
The boundary trace $\tr\big(\rho(K)\big)$ is:
\begin{equation*}
\kappa(x,y,z) = x^2 + y^2 + z^2 - x y z -2
\end{equation*}
so the relative character varieties are the level sets $\kappa^{-1}(t)$.
This mapping class group $\mcg$ acts by polynomial transformations of $\C^3$,
preserving the function $\kappa$ (compare \cite{puncturedtorus}). 
The Poisson structure is given by the bivector field
\begin{align*}
d\kappa \cdot \big(\dd{x}\wedge\dd{y}\wedge\dd{z}\big) & =
\big( 2x - y z\big) \dd{y}\wedge\dd{z} \\ & + 
\big( 2y - z x\big) \dd{z}\wedge\dd{x} \\ & + 
\big( 2z - x y\big) \dd{x}\wedge\dd{y}
\end{align*}
(where $\partial_x$ denotes $\frac{\partial}{\partial x}$, etc.).

\subsubsection{The four-holed sphere}
When $\Sigma$ is a four-holed sphere, the relative character varieties admit
a similar description. Present the fundamental group as
\begin{equation*}
\pi = \langle A, B, C, D  \mid  A B C D = 1 \rangle  
\end{equation*}
where the generators $A,B,C,D$ correspond to the components of $\partial\Sigma$.
The elements 
\begin{equation*}
X = A B, Y = B C, Z = C A
\end{equation*}
correspond to simple closed curves on $\Sigma$.
Denoting the trace functions $\hpgg \longrightarrow \C$ corresponding to elements
$A,B,C,D,X,Y,Z\in\pi$ by lower-case, the relative character varieties are
defined by:
\begin{align*}
x^2 + y^2 + z^2 + x y z & = (a b + c d) x  + (b c + a d) y \\
& + (a c + b d) z + (4 - a^2 - b^2  - c^2 - d^2  - a b c d)   
\end{align*}
with Poisson structure
\begin{align*}
\xi & =  
\big(a b + c d - 2x - y z\big) \dd{y}\wedge\dd{z} \\ & + 
\big(b c + d a - 2y - z x\big) \dd{z}\wedge\dd{x} \\ & + 
\big(c a + b d - 2z - x y\big) \dd{x}\wedge\dd{y}.
\end{align*}

\subsubsection{The two-holed torus}

Presenting the fundamental group of a two-holed torus as
\begin{align*}
\pi = \langle A, B, X, Y \mid  A X Y = Y X B \rangle, 
\end{align*}
where $A,B\in\pi$ correspond to the two components of $\partial\Sigma$, 
the elements 
\begin{align*}
Z &:= Y^{-1}X^{-1},\\
U &:= A X Y = B Y X, \\
V &: = B Y, \\
W &: = A X
\end{align*}
are represented by simple closed curves. Using the same notation for trace coordinates
as above, the relative character varieties are defined by the equations:
\begin{align*}
a + b & = x w + y v + u z - x y u \\
a b & = x^2 + y^2 + z^2 + u^2 + v^2 + w^2  \\
& \qquad  +  v w z - x y z - x u v - y u w - 4
\end{align*}
and the Poisson structure is 
\begin{align*}
(2 z - x y) \dd{x}\wedge\dd{y}  & + 
(2 x - y z) \dd{y}\wedge\dd{z}   + 
(2 y - y x) \dd{z}\wedge\dd{x}  + \\ 
(2 u - v x) \dd{v}\wedge\dd{x}  & + 
(2 v - x u) \dd{x}\wedge\dd{u}  + 
(2 x - u v) \dd{u}\wedge\dd{v}  + \\
(2 u - w y) \dd{w}\wedge\dd{y}  & + 
(2 w - y u) \dd{y}\wedge\dd{u}  + 
(2 y - u w) \dd{u}\wedge\dd{w}  +  \\ 
(2 (x y - z) & - v w)  \dd{v}\wedge\dd{w}      + 
(2 (x u - v) - w z)  \dd{w}\wedge\dd{z}  \\ & 
\qquad +  
(2 (y u - w) - z v) \dd{z}\wedge\dd{v}.
\end{align*}
These formulas are derived by applying the formulas for the Poisson 
bracket of trace functions developed in Goldman~\cite{invariantfunctions}
in combination with the trace identities in $\sltc$ (see \cite{fricke}).

\section{Compact Groups}

The simplest case occurs when $G=\Uo$. Then
\begin{equation*}
\hpgg = \hpg \cong \Uo^{2g} \cong H^1(\Sigma;\R/\Z)
\end{equation*}
is a $2g$-dimensional torus. If $M$ is a closed Riemann surface diffeomorphic
to $\Sigma$, then $\hpgg$ identifies with the {\em Jacobi variety\/} of $M$,
parametrizing topologically trivial holomorphic line bundles over $M$.
Although the complex structures on $\hpgg$ vary with the complex structures
on $\Sigma$, the symplectic structure
is independent of $M$. 


\subsection{Ergodicity}
The mapping class group action in this case factors through the
symplectic representation
\begin{equation*}
\mcg \longrightarrow \mathsf{Sp}(2g,\Z) 
\end{equation*}
(since the representation variety is just the ordinary cohomology
group with values in $\R/\Z$), which is easily seen to be ergodic.
This generalizes to arbitrary compact groups:

\begin{thm}
Let $G$ be a compact group. The $\Ou$-action on $\hpgg$ is ergodic.
\end{thm}
When the simple factors of $G$ are locally isomorphic to $\sut$ and
$\Sigma$ is orientable, this was proved in Goldman~\cite{erg}.  For
general $G$, this theorem is due to Pickrell-Xia~\cite{PickrellXia1}
when $\Sigma$ is closed and orientable, and
Pickrell-Xia~\cite{PickrellXia2} for compact orientable surfaces of
positive genus.

The following conjecture generalizes the above ergodicity phenomenon:
\begin{conjecture}\label{conj:forms}
Let $\Omega^*(\hpgg)$ 
be the de Rham algebra consisting of all 
measurable differential forms on $\hpgg$. Then the 
symplectic structures $\omega_\B$ generate the 
subalgebra of $\Omega^*(\hpgg)$ consisting of 
$\mcg$-invariant forms. 
\end{conjecture} 
\noindent
Since the $\mu$-measure of $\hpgg$ is finite, the representation
of $\mcg$ on 
\begin{equation*}
\HH :=  L^2(\hpgg,\mu))
\end{equation*}
is unitary. Andersen has informed me that he has
proved vanishing of the first cohomology group
$H^1(\mcg,\HH)$, and has raised the following conjecture
generalizing Conjecture~\ref{conj:forms}::

\begin{conjecture}
Suppose 
\begin{equation*}
C^\infty(\hpgg)\xrightarrow{D} C^\infty(\hpgg) 
\end{equation*}
is a differential
operator which commutes with the $\mcg$-action on $\hpgg$. Then
$D$ is a scalar multiple of the identity operator.
\end{conjecture}

\subsection{The unitary representation}
Ergodicity means that the only trivial subrepresentation of 
$\HH$ is the subspace $\C$ consisting of constants.
Furthermore the action is {\em weak mixing, \/}
by which we mean that $\C$ is the only {\em finite-dimensional \/}
invariant subspace~\cite{erg}. 
On the other hand the orthogonal complement $\HH_0$ to $\C$
in $\HH$ contains invariant subspaces.
For example the closure of the span of trace functions of 
{\em nonseparating simple closed curves on $\Sigma$\/} is 
an invariant subspace~\cite{reducible}.
\begin{prob} Decompose the representation on $\HH_0$ 
into irreducible representations of $\mcg$.
\end{prob}
When $G=\Uo$, and $\Sigma$ is the 2-torus, $\hpgg$ naturally
identifies with $T^2$, by the functions $\alpha,\beta$ corresponding
to a basis of $\pi_1(\Sigma)$. The functions 
\begin{equation*}
\phi_{m,n} := \alpha^m \beta^n,
\end{equation*}
forms a Hilbert basis of $\HH$, indexed by $(m,n)\in\Z^2$.
The $\mcg$-representation on $\HH$
arises from the linear $\gltz$-action on its basis $\Z^2$.
The $\gltz$-orbits on $\Z^2$ are indexed by integers $d\ge 0$.
The orbit of $(d,0)$ consists of all $(m,n)\in\Z^2$ with 
$\mathsf{gcd}(m,n) = d$. These are Hilbert bases for irreducible
constituents $C_d$ of $\HH$.

The irreducible constituents $C_d$ admit an alternate 
description, as follows.
The $d$-fold covering homomorphism 
\begin{equation*}
G\xrightarrow{\Phi_d} G  
\end{equation*}
induces a covering space
\begin{equation*}
\hpgg\longrightarrow \hpgg.  
\end{equation*}
Let $L_d$ denote the closure of the
image of the induced map $\HH \longrightarrow  \HH$. Then 
\begin{equation*}
L_d \,=\, \widehat{\bigoplus_{d'|d}} \; C_{d'} 
\end{equation*}
so $C_d$ consists of the orthocomplement in $L_d$ of the sum of all $L_{d'}$
for $d'|d$ but $d'\neq d$.

\begin{prob} 
Find a similar geometric interpretation
for the irreducible constituents for compact
nonabelian groups $G$.
\end{prob} 

\subsection{Holomorphic objects}
By Narasimhan-Seshadri~\cite{NarasimhanSeshadri}, and 
Ramanathan~\cite{Ramanathan},
a marked conformal structure $(f,M)$ on $\Sigma$
interprets $\hpgg$ as a moduli space of {\em holomorphic objects\/}
on $M$. To simplify the exposition we only consider the case $G=\Un$,
for which $\hpgg$ identifies with the moduli space $\mathfrak{U}_n(M)$ 
of {\em semistable  holomorphic $\C^n$-bundles \/}
over $M$ of zero degree~\cite{NarasimhanSeshadri}. 
The union of all $\mathfrak{U}_n(M)$ 
over $\langle f,M\rangle$ in $\Teich$ forms a holomorphic
fiber bundle
\begin{equation*}
\mathfrak{U}_n \longrightarrow \Teich
\end{equation*}
with an action of $\mcg$.
The quotient $\mathfrak{U}_n/\mcg$ fibers holomorphicly over $\Mod$. 
The Narasimhan-Seshadri theorem gives a (non-holomorphic) map
\begin{equation*}
\UU_n \xrightarrow{\hol} \hpgg
\end{equation*}
which on the fiber $\UU_n(M)$ is the bijection associating to
an equivalence class of semistable bundles the equivalence class of the
holonomy representation of the corresponding flat unitary structure.
$\UU_n/\mcg$ inherits a foliation $\FF_{\UU}$ from the
the  foliation of $\UU_n$ by level sets of $\hol$.
The dynamics of this foliation are equivalent to the dynamics of the
$\mcg$-action on $\hpg$.

Go one step further and replace $\Teich$ by its unit sphere bundle
$U\Teich$ and $\Mod$ by its (orbifold) unit sphere bundle
\begin{equation*}
U\Mod = (U\Teich)/\mcg. 
\end{equation*}
Pull back the fibration $\UU^k(\Sigma)$ to
$U\Mod$, to obtain a flat $\hpgg$-bundle $U\UU_n$
over $U\Mod$, 

The {\em Teichm\"uller geodesic flow\/} is a vector field on $U\Mod$
generating the geodesics for the Teichm\"uller metric on $\Teich$.
(Masur~cite{Masur}) Its horizontal lift with respect to the flat
connection is an vector field on the total space whose dynamics
mirrors the dynamics of the $\mcg$-action on $\hpgg$. 

As the $\mcg$-action on $\hpgg$
is weak-mixing, the unitary representation 
on $L^2\big(\hpgg,\mu \big)$ provides no
nontrivial finite-dimensional representations. Thus these
representations markedly differ from the representations
obtained by Hitchin~\cite{Hitchin2} and Axelrod, Della-Pietra, and
Witten~\cite{ADW} obtained from
projectively flat connections on $\hpgg$. 
Recently Andersen~\cite{Andersen} has proved that these finite-dimensional
projective representations of $\mcg$ are {\em asymptotically\/} faithful.

\subsection{Automorphisms of free groups}
Analogous questions arise for the outer automorphism group
of a free group $\pi$ of rank $r$. 
Let $G$ be a compact connected Lie group. 
Then Haar measure on $G$ defines 
an $\Out(\pi)$-invariant probability measure on $\hpg$.
\begin{conjecture}
If $r \ge 3$, the action of $\Out(\pi)$ on $\hpg$ is ergodic.
\end{conjecture}
\noindent
Using calculations in \cite{erg}, this conjecture has been proved~\cite{outfn}
when all of the simple factors of $G$ are locally isomorphic to $\sut$.

\subsection{Topological dynamics}
The topological theory is more subtle, since no longer may we ignore
invariant subsets of measure zero. For example, if $F\subset G$ is a finite
subgroup, then $\Hom(\pi,F)$ is finite and its image in 
$\hpgg$ is an invariant closed subset.

One might expect that if a representation $\rho\in\hpg$ has dense
image in $\sut$, that the $\mcg$-orbit of $[\rho]$ is dense
in $\hpgg$. This is true if $\Sigma$ is a one-holed torus
(Previte-Xia~\cite{PreviteXia1}) and 
if the genus of $\Sigma$ is positive
(Previte-Xia~\cite{PreviteXia2}).
In genus $0$, representations $\rho$ exist with dense
image but $\mcg\cdot [\rho]$ consists of only two points. 

Similar examples exist when $\Sigma$ is a four-holed sphere.
Benedetto and I showed~\cite{BG}, 
that when $-2 < a,b,c,d < 2$, the set of
$\R$-points of the relative character variety has one compact
component.
This component is diffeomorphic to
$S^2$. Depending on the boundary traces $(a,b,c,d)$, this component
corresponds to either $\sltr$-representations or
$\sut$-representations.  Previte and Xia~\cite{PreviteXia3} 
found representations $\rho$
in the components corresponding to $\sltr$-representations having
dense image, but whose orbit $\big(\mcg\cdot [\rho])$ has two points.
On the other hand, in both cases, Previte and Xia~\cite{PreviteXia4}
showed the action is {\em minimal\/} (every orbit is dense) for a
dense set of boundary traces in $[-2,2]^4$.

\begin{prob}
Determine necessary and sufficient conditions on a general representation
$\rho$ for its orbit $\mcg\cdot[\rho]$ to be dense.
\end{prob}
\noindent
The case when $G=\sut$ and $\Sigma$ an $n$-holed sphere for $n>4$ remains open.

\subsection{Individual elements}

For a closed surface of genus one, an individual element is ergodic
on the $\sut$-character variety if and only if it is {\em hyperbolic.\/}
In his doctoral thesis~\cite{RichBrown1,RichBrown2},
Brown used KAM-theory to show this
no longer holds for actions on relative $\sut$-character varieties over the
one-holed torus. Combining Brown's examples with a branched-cover 
construction suggests: 
\begin{prob}
Construct an example of a pseudo-Anosov mapping class for a closed surface
which is not ergodic on the $\sut$-character variety.
\end{prob}

\section{Noncompact Groups and Uniformizations}
For noncompact $G$, one expects less chaotic dynamics.
Trivial dynamics  -- in the form of {\em  proper\/} 
$\mcg$-actions --- occur for many
invariant open subsets corresponding to {\em locally homogeneous
geometric structures,\/}  (in the sense of Ehresmann~\cite{Ehresmann}) 
or {\em uniformizations.\/}  Such structures
are defined by local coordinate charts into a homogeneous space $G/H$ 
with coordinate changes which are restrictions of transformations from $G$.
Such an atlas globalizes  to a {\em developing map,\/} 
an immersion $\tilde\Sigma\longrightarrow G/H$
of the universal covering space $\tilde \Sigma \longrightarrow \Sigma$
which is equivariant with respect to a homomorphism $\pi\xrightarrow{\rho} G$.

To obtain a deformation space of such structures with an action of the
mapping class group, one introduces markings for a fixed topological
surface $\Sigma$, just as in the definition of Teichm\"uller space.
The {\em deformation space \/} $\Def$ consists of equivalence classes
of marked $(G,G/H)$-structures with a {\em holonomy map\/}
\begin{equation*}
\Def \xrightarrow{\hol} \hpgg 
\end{equation*}
which is $\mcg$-equivariant. The {\em Ehresmann-Thurston theorem\/}
asserts that, with respect to an appropriate topology on $\Def$, the
mapping $\hol$ is a local homeomorphism. (This theorem is implicit in
Ehresmann~\cite{Ehresmann2} and first explicitly stated by
Thurston~\cite{Thurston}. More detailed proofs were given by Lok~\cite{Lok},
Canary-Epstein-Green~\cite{CanaryEpsteinGreen}, and Goldman~\cite{geost}.
Bergeron and Gelander~\cite{BergeronGelander} 
give a detailed modern proof with applications to discrete subgroups.)

If $G=\pgltr$ and $G/H=\Ht$ is the hyperbolic plane, then $\Def = \Teich$.

Examples of uniformizations with proper $\mcg$-actions include:
\begin{itemize}
\item  $G=\psltr$:
The Teichm\"uller space $\Teich$, regarded as the component
of discrete embeddings in $\hpgg$;
\item 
$G=\psltc$:
Quasi-fuchsian space $\qf$ is an open subset of $\hpgg$ which 
is equivariantly biholomorphic to $\Teich\times\overline{\Teich}$;
\item 
$G = \slthr$.
The deformation space $\cc$ of convex $\rpt$-struc\-tures is a
connected component of $\hpgg$ (Choi-Gold\-man~\cite{ChoiGoldman1}) and
the $\mcg$-action is proper. More generally if $G$ is a split $\R$-form
of a semisimple group, Labourie~\cite{Labourie2}
has shown that $\mcg$ acts properly on the contractible component of $\hpgg$
discovered by Hitchin~\cite{Hitchin3}.
\end{itemize}

\subsection{Fricke-Teichm\"uller space}
A {\em Fuchsian representation\/} of $\pi$ into $G =\psltr$
is an isomorphism $\rho$ of $\pi=\pi_1(\Sigma)$ onto a discrete subgroup
of $G$. Since $\pi$ is torsionfree and $\rho$ is injective, $\rho(\pi)$ 
is torsionfree. Hence it acts freely on $\Ht$ and
the quotient $\Ht/\rho(\pi)$ is a complete hyperbolic surface. 
The representation $\rho$ defines a homotopy equivalence
\begin{equation*}
\Sigma \longrightarrow \Ht/\rho(\pi)
\end{equation*}
which is homotopic to a homeomorphism. Thus $\rho$ is the holonomy homomorphism
of a hyperbolic structure on $\Sigma$. 
The collection of $\pgltr$-conjugacy classes
of such homomorphisms identifies (via the Uniformization Theorem) 
with the {\em Teichm\"uller space\/} $\Teich$ of $\Sigma$.
When $\partial\Sigma\neq\emptyset$, then the {\em Fricke space\/} is
defined as the deformation space of complete hyperbolic structures on 
$\mathsf{Int}(\Sigma)$ such that each end is either a cusp or a complete
collar on a simple closed geodesic (a {\em funnel\/}). These representations
map each component of $\partial\Sigma$ to either a parabolic or a hyperbolic
element of $\psltr$ respectively. For details on Fricke spaces see
Bers-Gardiner~\cite{BersGardiner}.


The $\mcg$-action on $\Teich$ is proper.
This fact seems to have first been noted by Fricke~\cite{FrickeKlein}
(see Bers-Gardiner~\cite{BersGardiner} or Farb-Margalit~\cite{FarbMargalit}). 
It follows from two facts:
\begin{itemize}
\item $\mcg$ preserve a metric on $\Teich$;
\item The {\em simple marked length spectrum\/}
\begin{equation}\label{eq:markedlengthspectrum}
\bigg\{ \text{simple closed curves on~} \Sigma \bigg\}\bigg/ \Diffo  \;
\longrightarrow\; \R_+
\end{equation}
is a proper map.
\end{itemize}
See Abikoff~\cite{Abikoff}, Bers-Gardiner~\cite{BersGardiner},
Farb-Margalit~\cite{FarbMargalit} or
Harvey~\cite{Harvey} for a proof.  Another proof follows
from Earle-Eels~\cite{EarleEels}, and the closely related fact (proved
by Palais and Ebin~\cite{Ebin}) that the full diffeomorphism group of a
compact smooth manifold acts properly on the space of Riemannian
metrics. Compare \cite{GoldmanWentworth}.

\subsection{Other components and the Euler class}
Consider the case $G=\psltr$. 
Then the components of $\hpgg$ are indexed by the 
{\em  Euler class\/}
\begin{equation*}
\hpgg \stackrel{e}\longrightarrow H^2(\Sigma;\Z) \cong \Z
\end{equation*}
whose image equals
\begin{equation*}
\Z \;\cap\; [2-2g-b,2g-2+b]
\end{equation*}
where $\Sigma$ has genus $g$ and $b$ boundary components.
Thus $\hpgg$ has $4g+2b-3$ connected components
(\cite{Topcomps} and Hitchin~\cite{Hitchin1} when $b=0$). 
The main result of~\cite{doctoralthesis} is
that the two extreme components $e^{-1}\big(\pm(2-2g-b)\big)$
consist of discrete embeddings.
These two components differ by the choice of orientation, 
each one corresponding to $\Teich$, 
upon which $\mcg$ acts properly.
In contrast,
\begin{conjecture}\label{conj:ergr}
Suppose that $b=0$ ($\Sigma$ is closed).
For each integer $1 \le k \le 2g+b-2$, 
the $\mcg$-action on the component $e^{-1}(2-2g + b+ k)$ of $\hpg$
is ergodic.
\end{conjecture}
\noindent
When $b=0$, the component 
\begin{equation*}
e^{-1}(3-2g) \;\approx\; \Sigma\,\times\,\R^{6g-8}  
\end{equation*}
represents a $6g -6$-dimensional
{\em thickening\/} of $\Sigma$, upon which $\mcg$ acts. However,
Morita~\cite{Morita} showed that $\mcg$ cannot act smoothly on 
$\Sigma$ itself inducing the homomorphism $\Diff(\Sigma)\longrightarrow\mcg$.
(Recently Markovic~\cite{Markovic} has announced that if $\Sigma$ is 
a closed surface of genus $>5$, then 
$\mcg$ cannot even act on $\Sigma$ by homeomorphisms
inducing $\Homeo(\Sigma)\longrightarrow\mcg$.)

\begin{prob} 
Determine the smallest dimensional manifold homotopy-equivalent to $\Sigma$
upon which $\mcg$ acts compatibly with the outer action of $\mcg$ on
$\pi_1(\Sigma)$.
\end{prob} 

\subsection{The one-holed torus}\label{sec:oneholedtorus}

For surfaces with nonempty boundary, the dynamics appears more
complicated.  When $\Sigma$ is a one-holed torus ($g = b =1$) 
and $G=\psltr$ or $\sut$, 
this was completely analyzed in \cite{puncturedtorus}.

As in \S\ref{subsec:oneholedtorus}, 
the $\sltc$-character variety identifies with $\C^3$, where the three
coordinates $(x,y,z)$ are traces of three generators of $\pi$ corresponding
to the generators $X,Y,XY$.
In these coordinates, 
the trace of the element $K = XYX^{-1}Y$ of $\pi$ corresponding to 
$\partial\Sigma$ equals
\begin{equation*}
\kappa(x,y,z) := x^2 + y^2 + z^2 - x y z - 2.
\end{equation*}
The {\em relative $\sltc$-character variety\/} of $\Sigma$ is then
the family of level sets $\kappa^{-1}(t)$ of $\C^3\xrightarrow{\kappa}\C$. 

The set $\kappa^{-1}(t)\cap \R^3$ of $\R$-points of $\kappa^{-1}(t)$, 
for boundary trace $t\in\R$, are of two types:
\begin{itemize}
\item 
The $\sut$-characters, with $x,y,z\in[-2,2]$ and  $t < 2$;
\item 
The $\sltr$-characters, with either:
\begin{itemize} 
\item at least one of  $x,y,z$ lies in $(-\infty,-2]\,\cup\,[2,\infty)$,  or
\item each $x,y,z$ lies in $[-2,2]$ and $t\ge 2$.
\end{itemize}
\end{itemize}
If $\vert t\vert > 2$, no $\sut$-characters lie in 
$\kappa^{-1}(t)\cap\R^3$.
If $t\neq 2$, these two subsets of $\kappa^{-1}(t)\cap\R^3$ 
are disjoint. If $t=2$, these two subsets intersect on the subset 
\begin{equation*}
[-2,2]^3 \cap \kappa^{-1}(2) 
\end{equation*}
corresponding to $\sot$-representations. The space
of $\sot$-characters is 2-sphere with 4 branch points of cone angle $\pi$  
(a tetrahedron with smoothed edges).

The $\mcg$-action determines a dynamical system on each level set.
By Keen~\cite{Keen},
the Fricke space of $\Sigma$ is the subset 
\begin{equation*}
\big\{ (x,y,z)\in\R^3 \mid \kappa(x,y,z) \le -2\big\}
\end{equation*}
with a proper $\mcg$-action.
Each level set $\R^3\cap\kappa^{-1}(t)$, for $t < -2$, 
is homeomorphic to a disjoint union of four discs; 
the four components are distinguished by different lifts of the
representation from $\psltr$ to $\sltr$. 

The level set $\R^3\cap\kappa^{-1}(-2)$ has one notable feature.
It has five components, four of which correspond
to the Teichm\"uller space of $\Sigma$, and the other component 
$\{(0,0,0)\}$ consists of just the origin. The Teichm\"uller space
(corresponding to the deformation space of complete {\em finite area\/}
hyperbolic structures) corresponds to representations taking the
boundary element of $\pi$ to a parabolic transformation of trace $-2$.
On the other hand, $\{(0,0,0)\}$ corresponds the 
{\em quaternion representation\/} in $\sut$: 
\begin{align*}
X & \longmapsto  \bmatrix i & 0 \\ 0 & -i\endbmatrix,  \\
Y & \longmapsto  \bmatrix 0 & -1 \\ 1 & 0\endbmatrix 
\end{align*}
The peripheral element $K\in\pi$ 
maps to the nontrivial {\em central\/} element $-I\in\sut$. 

Here we see --- for the first time --- 
the coexistence of two extremes of dynamical behavior: 
\begin{itemize}
\item The proper action on the $\sltr$-characters;
\item The entire mapping class group $\mcg$ fixes a point,
in a sense, the ``most chaotic'' action.
\end{itemize}

This dichotomy persists for $-2 < t < 2$. 
The origin deforms to a
compact component, consisting of characters of $\sut$-representations
with an ergodic  $\mcg$-action. Four contractible components,
correspond to holonomy representations of hyperbolic structures
on a torus with a cone point. The cone angle $\theta$ relates to 
the boundary trace by
\begin{equation*}
t = -2 \cos(\theta/2). 
\end{equation*}
The $\mcg$-action on these components is proper.

Although $\mcg$ acts properly, none of the corresponding
representations are discrete embeddings. The key property seems to
be that nonseparating simple loops are mapped to hyperbolic elements,
so the simple marked length spectrum \eqref{eq:markedlengthspectrum}
is a proper map.

\begin{prob} Find general conditions which ensure that 
\eqref{eq:markedlengthspectrum} is proper.
\end{prob}

The level set $\R^3\cap\kappa^{-1}(2)$ consists of characters of
abelian representations, and $\mcg$ is ergodic on each of the four
connected components of the smooth part of $\R^3\cap\kappa^{-1}(2)$.
When $2 < t \le 18$, the $\mcg$-action on $\R^3\cap\kappa^{-1}(t)$ is
ergodic.

For $t > 18$, the level sets $\R^3\cap\kappa^{-1}(t)$ display both
proper dynamics and chaotic dynamics. The region $(-\infty,-2]^3$
consists of characters of discrete embeddings $\rho$ where the
quotient hyperbolic surface $\Ht/\rho(\pi)$ is homeomorphic to a
three-holed sphere.  Every homotopy equivalence $\Sigma\longrightarrow
P$, where $P$ is a hyperbolic surface homeomorphic to a three-holed
sphere, determines such a character. Furthermore these determine closed
triangular regions which are freely permuted by $\mcg$. On the complement
of these wandering domains the action is ergodic.

When $G=\pgltr$, the group of (possibly orientation-reversing) isometries
of $\Ht$, a similar analysis was begun by 
Stantchev~\cite{Stantchev,GoldmanStantchev}. 
One obtains similar dynamical systems,
where $\mcg$ acts now on the space of representations into the group 
\begin{equation*}
G_\pm \;=\; \sltc \,\cap\, \big(\gltr\, \cup\, i\ \gltr\big) 
\end{equation*}
which doubly covers the two-component group $\pgltr$.  
These $G_\pm$-representations are again
parametrized by traces. They comprise four components, 
one of which is the subset of $\R^3$ parametrizing $\sltr$-representations
discussed above. The other three components are 
\begin{equation*}
\R\times i\R\times i\R,\quad i\R\times \R\times i\R,\quad 
i\R\times i\R\times \R 
\end{equation*}
respectively.  Consider 
$i\R\times \R\times i\R$. 
For $-14 \le t < 2$, the $\mcg$-action is ergodic, but
when $t < -14$, wandering domains appear. 
The wandering domains  correspond to
homotopy-equivalences $\Sigma\longrightarrow P$, where $P$ is a
hyperbolic surface homeomorphic to a two-holed projective plane.
The action is ergodic on the complement of the wandering domains.

\begin{prob} Determine the ergodic behavior of the $\mcg$-action
on the level sets 
\begin{equation*}
\big(i\R\times \R\times i\R \big) \cap \kappa^{-1}(t)
\end{equation*}
where $t > 2$. The level sets for $t>6$ contains wandering domains
corresponding to Fricke spaces of a one-holed Klein bottle.
\end{prob}

\subsection{Hyperbolic 3-manifolds}

When $G = \psltc$, the subset $\qf$ of $\hpgg$ corresponding to 
embeddings of $\pi$ onto quasi-Fuchsian subgroups of $G$ is
open and $\mcg$-invariant. Furthermore the Bers isomorphism~\cite{Bers} 
provides a $\mcg$-invariant biholomorphism
\begin{equation*}
\qf \longrightarrow \Teich\times \overline{\Teich}. 
\end{equation*}
Properness of the action of $\mcg$  on $\Teich$ implies properness
on $\qf$.
 
Points on the boundary of $\qf$ also correspond to discrete
embeddings, but the action is much more complicated. Recently Souto
and Storm~\cite{SoutoStorm} have proved that $\partial\qf$ contains
a $\mcg$-invariant closed nowhere dense topologically perfect set
upon which the action is topologically transitive. From this they
deduce that every continuous $\mcg$-invariant function on $\partial\qf$
is constant.

While for representations into $G=\psltr$, the $\mcg$-orbits of
discrete embeddings are themselves {\em discrete,\/} the situation
becomes considerably more complicated for larger $G$. For $G=\psltc$,
representations corresponding to the fiber of a hyperbolic mapping
torus furnish points with infinite stabilizer.  This is one of the
easiest ways to see that $\mcg$ does not act properly on characters
of discrete embeddings.  Namely, if $M^3$ is a
hyperbolic 3-manifold which admits a fibration $M^3 \xrightarrow{f} S^1$,
then the class of the restriction $\rho$ of the holonomy
representation
\begin{equation*}
\pi_1(M^3)\longrightarrow \psltc 
\end{equation*}
to the surface group
\begin{equation*}
\pi \;:=\; \pi_1(f^{-1}(s_0)) \;\cong\; 
\Ker\big(\pi_1(M) \xrightarrow{f_*} \Z\big) 
\end{equation*}
is invariant under the monodromy automorphism $h\in\Aut(\pi)$ of
$M^3$. That is, there exists $g\in\psltc$ such that 
\begin{equation*}
\rho\big(h(\gamma)\big)  = g\rho(\gamma) g^{-1}
\end{equation*}
for all $\gamma\in\pi$.
Furthermore $[\rho]$ is a smooth point of $\hpgg$.
Kapovich~\cite{Kapovich} proved McMullen's conjecture~\cite{McMullen1} 
that the  derivative of the mapping class $[h]$ at $[\rho]$ is hyperbolic,
that is, no eigenvalue has norm 1. This contrasts the case of
{\em abelian representations,\/} since homologically trivial 
pseudo-Anosov mapping classes act trivially on $\hpg$ .

Thus $\mcg$ does not act properly 
on the set of characters of discrete embeddings. 
Let $[\rho]$ (as above) be the character of a discrete embedding of $\pi$ 
as the fiber of a hyperbolic mapping torus.
The stabilizer of $[\rho]$ contains the infinite cyclic group generated by 
the mapping class corresponding to $[h]$. (In fact $\langle[h]\rangle$ 
has finite index in the stabilizer of $[\rho]$.) Since stabilizers of
proper actions of discrete groups are finite, $\mcg$ does not act
properly.

The Souto-Storm theorem shows that this chaotic dynamical behavior
pervades the entire boundary of quasi-Fuchsian space $\qf$.

In another direction, using ideas generalizing those of
Bowditch~\cite{Bowditch},  Tan, Wong and Zhang~\cite{TWZ} have shown
that the action of $\mcg$ on the representations satisfying the analogue
of {\em Bowditch's Q-conditions\/} is proper. This also generalizes
the properness of the action on the space of quasi-Fuchsian
representations.

At present little is known about the dynamics of $\mcg$ acting
on the $\sltc$-character variety. Conversations with Dumas led to
the following problem:

\begin{prob}
Find a point $\rho\in\Hom(\pi,\sltc)$ such that the closure
of its orbit $\overline{\mcg\cdot[\rho]}$ meets both 
the image of the unitary characters $\Hom(\pi,\sut)$ 
and the closure $\overline{\qf}$ of the quasi-Fuchsian characters.
\end{prob}

\subsubsection{Homological actions}
The action of $\mcg$ on the homology of $\hpgg$ furnishes another
source of possibly interesting linear representations of $\mcg$.
With Neumann~\cite{goldmanneumann}, we proved that for the
relative $\sltc$-character varieties of the one-holed torus and
four-holed sphere, the action of $\mcg$ factors through a finite group. 

Atiyah-Bott~\cite{AtiyahBott} use infinite-dimensional Morse theory to analyze
the algebraic topology of $\hpgg$, when $G$ is compact. 
For the nonsingular components their techniques imply that 
the $\mcg$-action on the rational cohomology of $\hpgg$  factors through
the symplectic representation of $\mcg$ on $H^*(\Sigma)$.
In particular Biswas~\cite{Biswas} proved that the 
Torelli group acts trivially on nonsingular components.
In contrast, Cappell-Lee-Miller~\cite{CLM1,CLM2} proved the surprising
result that that the Torelli group acts {\em nontrivially\/}
on the homology of the $\sut$-character variety when $\Sigma$ is closed.

\subsection{Convex Projective Structures and Hitchin representations}

When $G=\slthr$, the mapping class group $\mcg$ acts properly on the
component $\cc$ of $\Hom(\pi,G)/G$ corresponding to convex
$\rpt$-structures (Goldman~\cite{convex},
Choi-Goldman~\cite{ChoiGoldman1,ChoiGoldman2}).  This component is
homeomorphic to a cell of dimension $-8\chi(\Sigma)$, and, for a
marked Riemann surface $M$ homeomorphic to $\Sigma$, admits the
natural structure of a holomorphic vector bundle over $\Teich$. The
work of Labourie~\cite{Labourie1} and
Loftin~\cite{Loftin1,Loftin2,Loftin3} gives a more intrinsic
holomorphic structure on $\cc$.

The existence of this contractible component is a special case of a
general phenomenon discovered by Hitchin~\cite{Hitchin3}. Hitchin
finds, for {\em any\/} split real form of a semisimple group $G$, a
contractible component in $\hpgg$. For $G = \SL(n,\R)$ this
component is characterized as the component containing the composition
of discrete embeddings 
\begin{equation*}
\pi\longrightarrow\SL(2,\R) 
\end{equation*}
with the
irreducible representation 
\begin{equation*}
\SL(2,\R)\longrightarrow\SL(n,\R). 
\end{equation*}
Recently, Labourie has found a dynamical description~\cite{Labourie2}
of representations in Hitchin's component, and 
has proved they  are discrete embeddings.
Furthermore he has shown that $\mcg$ acts properly on this component. 
(These closely relate to the 
{\em higher Teichm\"uller spaces \/}
of Fock-Goncharov~\cite{FockGoncharov1,FockGoncharov2}.)

When $G$ is the automorphism group of a Hermitian symmetric space of
noncompact type, Bradlow, Garcia-Prada, and Gothen have found other
components of $\hpgg$, for which the {\em Toledo invariant,\/} is
maximal~\cite{BradlowGarciaPradaGothen1,BradlowGarciaPradaGothen2}.
Their techniques involve Morse theory along the lines of 
Hitchin~\cite{Hitchin1}. 
Recently, Burger, Iozzi, and Wienhard
have shown~\cite{BurgerIozziWienhard} that the representations
of maximal Toledo invariant consist of discrete embeddings.
Using results from \cite{BurgerIozziLabourieWienhard},
Wienhard has informed me that $\mcg$ acts properly on these components,

\subsection{The energy of harmonic maps}
An interesting invariant of surface group representations
arises from the theory of {\em twisted harmonic maps\/} of Riemann
surfaces, developed in detail in collaboration with
Wentworth~\cite{GoldmanWentworth}.
Namely to each
{\em reductive representation\/} $\pi\xrightarrow{\rho} G$,
one associates an {\em energy function\/}
\begin{equation*}
\Teich \xrightarrow{E_\rho} \R  
\end{equation*}
whose qualitative properties reflect the $\mcg$-action.
Assuming that the Zariski closure of $\rho(\pi)$ in $G$ is reductive,
for every marked Riemann surface $\Sigma \longrightarrow M$, there
is a $\rho$-equivariant harmonic map
\begin{equation*}
\tilde{M} \longrightarrow G/K 
\end{equation*}
where $K\subset G$ is a maximal compact subgroup
(Corlette~\cite{Corlette}, Donaldson~\cite{Donaldson}, 
Labourie~\cite{Labourie}, and Jost-Yau \cite{Jost-Yau},
following earlier work by Eels-Sampson~\cite{EelsSampson}).
Its energy density determines an exterior $2$-form on $\Sigma$, 
whose integral is defined as $E_\rho(\langle f,M\rangle)$.

When $\rho(\pi)$ lies in a compact subgroup of $G$, then the
twisted harmonic maps are constant, and the energy function is
constantly zero. At the other extreme is the following
result, proved in \cite{GoldmanWentworth}:

\begin{thm}\label{thm:gw}
Suppose that $\rho$ is an embedding of $\pi$ onto a convex cocompact discrete
subgroup of $G$. Then the energy function $E_\rho$ is a proper function
on $\Teich$.
\end{thm}
Here a discrete subgroup $\Gamma\subset G$ is {\em convex cocompact\/}
if there exists a geodesically convex subset $N\subset G/K$ such that
$\Gamma\backslash N$ is compact. For $\psltc$, these are just the
quasi-Fuchsian representations.  This result was first proved by
Tromba~\cite{Tromba} for Fuchsian representations in $\psltr$, and the
ideas go back to Sacks-Uhlenbeck~\cite{SacksUhlenbeck} and
Schoen-Yau~\cite{SchoenYau}. 

It is easy to prove (see \cite{GoldmanWentworth}) that if $\Omega\subset\hpgg$
is a $\mcg$-invariant open set for which each function $E_\rho$ is proper,
for $[\rho]\in\Omega$, then the action of $\mcg$ on $\Omega$ is proper.
This gives a general analytic condition implying properness.

Unfortunately, convex cocompactness is extremely restrictive; 
Kleiner and Leeb have proved~\cite{KleinerLeeb} that in rank $>1$
such groups are {\em never\/} Zariski dense. However, we know many 
examples (the deformation space $\cc$ of convex $\rpt$-structures,
the Hitchin representations by Labourie~\cite{Labourie2}, other
components of maximal representations~\cite{BradlowGarciaPradaGothen1,
BradlowGarciaPradaGothen2, BurgerIozziLabourieWienhard})
where we expect the $\mcg$-action to be proper.
The only use of geodesic convexity in the above result is that
the images of harmonic maps are constrained to lie in the set $N$.

\begin{prob}
Find a substitute for convex cocompactness in higher rank which includes
the above examples of proper $\mcg$-actions, and for which $E_\rho$ is
proper.
\end{prob}

The work of Bonahon-Thurston on geometric tameness, and its
recent extensions, implies that the energy function of
a discrete embedding $\pi\longrightarrow\psltc$
is proper if and only if it is quasi-Fuchsian~\cite{GoldmanWentworth}.

\subsection{Singular uniformizations and complex projective 
structures}

When $G=\psltr$, the other components of $\Hom(\pi,G)$ may be studied
in terms of hyperbolic structures with singularities as follows.
Instead of requiring all of the coordinate charts to be local
homeomorphisms, one allows charts which at isolated points look like
the map
\begin{align*}
\C & \longrightarrow \C  \\
z & \mapsto z^k \end{align*}
that is, the geometric structure has an isolated singularity of 
{\em cone angle\/}
$\theta=2k\pi$.
Such a singular hyperbolic structure may be alternatively described
as a singular Riemannian metric $g$ whose
curvature equals $-1$ plus Dirac distributions weighted by
$2\pi-\theta_i$ at each singular point $p_i$ of cone angle $\theta_i$.
The structure is nonsingular on the complement 
$\Sigma - \{p_1,\dots,p_k\}$, and that hyperbolic structure has
holonomy representation
\begin{equation*}
\pi_1\big(\Sigma - \{p_1,\dots,p_k\}\big) \xrightarrow{\hat{\rho}} \psltr
\end{equation*}
such that the holonomy of a loop $\gamma_i$ encircling $p_i$ is elliptic
with rotation angle $\theta_i$.

In particular if each $\theta_i\in 2\pi\Z$, then
$\hat\rho(\gamma_i)=1$. The representation $\hat\rho$ extends
to a representation $\rho$ of $\pi_1(\Sigma)$:
\begin{equation*}
\xymatrix{
& \pi_1\big(\Sigma - \{p_1,\dots,p_k\}\big) 
\ar@{>>}[d] \ar[dr]^{\hat\rho} \\
  & \pi_1(\Sigma) \ar@{-->}[r]_{\rho} & \psltr.           }
\end{equation*}
Applying Gauss-Bonnet to $g$ implies that $\rho$ has
Euler class
\begin{equation*}
e(\rho) = \chi(M)+\frac{1}{2\pi} \sum_{i=1}^{k} (\theta_i-2\pi).             
\end{equation*}

It is convenient to assume that each $\theta_i=4\pi$ and the points
$p_i$ are not necessarily distinct --- 
a cone point of cone angle $4\pi$ with multiplicity $m$ 
is then a cone point with cone angle $2(m+1)\pi$.
The uniformization theorem of McOwen~\cite{McOwen}, 
Troyanov~\cite{Troyanov}, and Hitchin~\cite{Hitchin1} asserts:
given a Riemann surface $M\approx \Sigma$,  
there exists a unique singular hyperbolic structure in the conformal class of
$M$ with cone angle $\theta_i$ at $x_i$ for $i=1,\dots,k$ as long as 
\begin{equation*}
\chi(\Sigma)+\frac{1}{2\pi} \sum_{i=1}^{k} (\theta_i-2\pi) < 0.  
\end{equation*}
(Hitchin only considers the case when $\theta_i$ are multiples of
$2\pi$, while McOwen and Troyanov deal with arbitrary positive angles.)
The resulting {\em uniformization map \/}  assigns to the
collection of points $\{p_1,\dots,p_k\}$ (where $0\le k\le |\chi(\Sigma)|$) 
the singular hyperbolic structure with cone angles $4\pi$ (counted
with multiplicity) at the $p_i$.
The equivalence class of the holonomy representation in the component
\begin{equation*}
e^{-1}(\chi(\Sigma)+k)\subset\hpgg
\end{equation*}
defines a map from the
symmetric product $\Sym^k(M)$ to $e^{-1}(\chi(M)+k)$.
The following result follows from Hitchin~[35]:

\begin{thm} Let $M$ be a closed Riemann surface. The above map 
\begin{equation*}
\Sym^k(M)\longrightarrow e^{-1}(\chi(M)+k)     
\end{equation*}
is a homotopy equivalence. \end{thm}

The union of the symmetric powers $\Sym^k(M)$, one for each marked
Riemann surface $M$, over $\langle f,M\rangle\in\Teich$, can be given
the structure of a holomorphic fiber bundle $\SSS^k(\Sigma)$ over $\Teich$,
to which the action of $\mcg$ on $\Teich$ lifts.
The above maps define a homotopy equivalence
\begin{equation*}
\SSS^k(\Sigma)\stackrel{\Uu}\longrightarrow e^{-1}(\chi(\Sigma)+k).     
\end{equation*}
which is evidently $\mcg$-equivariant. However, since $\mcg$ acts
properly on $\Teich$, 
it also acts properly on the $(6g-6+2k)$-dimensional
space $\SSS^k(\Sigma)$. The quotient
$\SSS^k(\Sigma)/\mcg$  is the total space of an 
(orbifold)  $\Sym^k(\Sigma)$-bundle  over the
Riemann moduli space
\begin{equation*}
\Mod := \Teich/\mcg.
\end{equation*}
The fibers of $\Uu$ define a (non-holomorphic)
foliation 
of $\SSS^k(\Sigma)/\mcg$, 
{\em a flat $\Sym^k(\Sigma)$-bundle,\/}
whose dynamics mirrors the dynamics of the $\mcg$-action on the 
component $e^{-1}(\chi(\Sigma) + k)$.

In general, $\Uu$ is {\em not\/} onto: if $\Sigma\xrightarrow{f}\Sigma'$
is a degree one map to a closed surface $\Sigma'$ of smaller genus,
and $\rho'$ is a Fuchsian representation of $\pi_1(\Sigma')$, then
$\rho :=\rho'\circ f_*$ lies outside $\Image(\Uu)$.
The following conjecture arose in discussions with Neumann:

\begin{conjecture}
If $k=1$, then $\Uu$ is onto. 
In general a $\psltr$-representation with dense image lies in 
$\Image(\Uu)$.
\end{conjecture}

\subsection{Complex projective structures}

A similar construction occurs with the deformation space $\cpo(\Sigma)$
of {\em marked $\cpo$-structures\/} on $\Sigma$. A $\cpo$-manifold
is a manifold with a coordinate atlas modeled on $\cpo$, with coordinate
changes in $G = \psltc$. The space $\cpo(\Sigma)$ consists of equivalence
classes of marked $\cpo$-structures, that is, homeomorphisms 
$\Sigma\longrightarrow N$ where $N$ is a $\cpo$-manifold. Since $\psltc$
acts holomorphicly, the $\cpo$-atlas is a holomorphic atlas and $N$
has a underlying Riemann surface $M$. The resulting $\mcg$-equivariant map
\begin{equation*}
\cpo(\Sigma) \longrightarrow \Teich 
\end{equation*}
is a holomorphic affine bundle, 
whose underlying vector bundle is
the holomorphic cotangent bundle of $\Teich$. 
In particular $\mcg$ acts properly on $\cpo(\Sigma)$ with quotient
a holomorphic affine bundle over $\Mod$.

The map which associates to a marked $\cpo$-structure on $\Sigma$ its
holonomy representation is a local biholomorphism
\begin{equation*}
\cpo(\Sigma) \xrightarrow{\hol} \ \hpgg
\end{equation*}
which is known to be very
complicated. Gallo-Kapovich-Marden~\cite{GalloKapovichMarden} have
shown that its image consists of all equivalence classes of representations 
$\rho$ for which:
\begin{itemize}
\item $\rho$ lifts to a representation $\pi\longrightarrow\sltc$;
\item The image $\rho(\pi)$ is not precompact;
\item The image $\rho(\pi)$ is not solvable.
\end{itemize}
The latter two conditions are equivalent to $\rho(\pi)$ not leaving
invariant a finite subset of $\Hthree \cup \partial\Hthree$. 
(The cardinality of this finite subset can be taken to be either $1$ or
$2$.)

The holonomy map $\hol$ is $\mcg$-equivariant. The action of 
$\mcg$ on $\cpo(\Sigma)$ is proper, since it covers the action of
$\mcg$ on $\Teich$. The quotient 
$\cpo(\Sigma)/\mcg$ affinely fibers over $\Mod$. As before $\hol$ defines
a foliation of $\cpo(\Sigma)/\mcg$ orbit equivalent to the 
$\mcg$-action on $\hpgg$. Thus $\hol$ may be regarded as a 
{\em resolution \/} of the $\mcg$-action.

As a simple example, the $\mcg$-action is proper on the quasi-Fuchsian
subset $\qf\subset\hpgg$. As noted above, it is a maximal open set
upon which $\mcg$ acts properly.
Its restriction
\begin{equation*}
\hol^{-1}(\qf) \longrightarrow \qf 
\end{equation*}
is a covering space (\cite{ProjSt}). 
However, the bumping phenomenon discovered by
McMullen~\cite{McMullen2} implies that $\hol$ is not a covering space
on any open neighborhood strictly containing $\qf$.


\end{document}